\documentclass[preprint, 10pt]{imsart}

\RequirePackage[OT1]{fontenc}
\RequirePackage{amsthm,amsmath}
\RequirePackage[authoryear,round]{natbib}
\RequirePackage[colorlinks,citecolor=blue,urlcolor=blue,pdfstartview=FitB]{hyperref}

\usepackage{tikz}
\usetikzlibrary{matrix}
\usepackage{multicol}
\usepackage{changepage}
\usepackage{float}
\usepackage{amsmath,amssymb,epsfig}
\usepackage{graphics}
\usepackage{tabularx}
\usepackage{framed}
\usepackage[mathscr]{euscript}
\usepackage{imsart}

\usepackage{epstopdf}

\usepackage{titlesec}
\titleformat*{\section}{\large\bfseries}

\usepackage{graphicx}
\usepackage{bbm}

\startlocaldefs
\numberwithin{equation}{section}
\theoremstyle{plain}

\long\def\comment#1{}

\theoremstyle{definition}

\newcommand{\eps}{\varepsilon}

\newcommand{\be}{\begin{eqnarray}}
\newcommand{\ee}{\end{eqnarray}}

\newcommand{\ba}{\begin{array}}
\newcommand{\ea}{\end{array}}
\newcommand{\bs}{\begin{align}\begin{split}\nonumber}
\newcommand{\bsnumber}{\begin{align}\begin{split}}
\newcommand{\es}{\end{split}\end{align}}

\renewcommand{\(}{\left(}
\renewcommand{\)}{\right)}

\renewcommand{\hat}{\widehat}

\newcommand{\QED}{$\blacksquare$}

\newcommand{\Ep}{{\mathrm{E}}}

\renewcommand{\Pr}{{\mathrm{P}}}

\renewcommand{\hat}{\widehat}
\renewcommand{\leq}{\leqslant}
\renewcommand{\geq}{\geqslant}

\newcommand{\kmax}{{k\text{-}\mathrm{max}}}
\newcommand{\ktildemax}{{k\text - \widetilde \max }}

\usepackage[flushleft]{threeparttable}
\usepackage{multirow}
\usepackage{algorithm,algorithmic}
\usepackage{mathrsfs}
\usepackage{enumitem}   
\usepackage{subfig}
\usepackage{pdflscape}
\endlocaldefs

\begin{document}

\begin{frontmatter}

{\large \textbf{Dimension-Free Anticoncentration Bounds for Gaussian Order Statistics with Discussion of Applications to Multiple Testing}}

\title{} 
\runtitle{Dimension-Free Anticoncentration Bounds for Gaussian Order Statistics}

\runauthor{D. Kozbur \ }
\author{\fnms{Damian} \snm{Kozbur}\ead[label=e3]{damian.kozbur@econ.uzh.ch}}
\thankstext{t1}{\today.}
\affiliation{University of Z\"urich }
\address{  Department of Economics, University of Z\"urich \\ Sch\"onberggasse 1, 
8001 Z\"urich, Switzerland\\
\printead{e3}\\
}

\

\begin{abstract} The following anticoncentration property is proved.  The probability that the $k$-order statistic of an arbitrarily correlated jointly Gaussian random vector $X$ with unit variance components lies within an interval of length $\eps$ is bounded above by $2{\eps}k ({ 1+\Ep[\|X\|_\infty ]}) $.  This bound has implications for generalized error rate control in statistical high-dimensional multiple hypothesis testing problems, which are discussed subsequently. 
\end{abstract}

\end{frontmatter}

\noindent \textbf{Introduction}  

\

\noindent Consider a random vector $X$, taking values in $\mathbb R^{p}$ for some positive integer $p$, with components denoted $X_1,...,X_p$.  Let $\kmax(X)$ denote the value of the $k$-th largest component of $X$.  The following theorem holds.

\

\noindent
\textrm{THEOREM 1.} \textit{Let $X \sim \mathsf N(0,\Sigma) \in \mathbb R^{p}$ be a Gaussian random vector with unit variance components, ie. the only restriction on $\Sigma$ is that $\Ep[X_1^2]=...=\Ep[X_p^2]=1$.  Then} $$\sup_{y \in \mathbb R} \Pr( k\text{-}\mathrm{max}(X) \in [y,y+\eps]) \leq {2\eps} k ({ 1+\Ep[\|X\|_\infty ]}) .$$ 

\

\noindent Theorem 1 is proved below after a literature review and general description of key new arguments needed for its proof.  Theorem 1 has implications for generalized error rate control in statistical high-dimensional multiple hypothesis testing problems, which are then discussed.

\

\noindent Bounds of the form of Theorem 1---upper bounds on the probability that a random variable falls in intervals of small width---are referred to commonly as anticoncentration bounds.  The study of anticoncentration bounds for sums of independent random variables dates back to at least \cite{littlewood:offord} and {\cite{erdos:anticoncentration}} with the later reference proving that $\sup_{y \in \mathbb R} \Pr (\sum_{i=1}^pz_j R_j \in [y,y+2] ) \leq 2^{-p} \hspace{.5mm} {_pC}_{\lfloor p/2 \rfloor}$ for any fixed real $z_j$ with $|z_j| \geq 1$ and  $R_j \in \{\pm 1\}$ independent Rademacher random variables.

\

\noindent In the case of $k=1$, the conclusion of Theorem 1 was proved in \cite{chern:chet:kato:anti}.   Note that the bound in Theorem 1 is dimension-free (ie. independent of $p$).  Thus it is potentially much tighter than dimension-dependent bounds when components $X_j$ are sufficiently correlated.  Prior to \cite{chern:chet:kato:anti}, \cite{nazarov:anticonc} proved dimension-dependent anticoncentration bounds in the case $k=1$.  \cite{anticoncentration:malliavin} give a dimension-free bound in the case $k=1$ but required that all components of $\Sigma$ be strictly positive.  

\

\noindent  The bound for $k=1$ has found various applications to statistics.  It is a major component in proving asymptotic family-wise error rate control in a variety of high-dimensional statistical hypothesis testing problems using the multiplier bootstrap as in \cite{chet:chern:kato:multboot} and using the empirical bootstrap as in \cite{clt:empboot:highdim}.  In the case of $k=1$, \cite{CHERNOZHUKOV20163632} and \cite{Chernozhukov_2014} used the dimension-free property in application to studying infinite dimensional Gaussian elements and  empirical process theory.  Anticoncentration bounds are also used extensively in the course of proving various high-dimensional central limit theorems in \cite{chernozhukov2019improved}, \cite{lopes2020central}, \cite{yuta:clt}, and \cite{chernozhukov2021nearly}.

\

\noindent The case of $k\geq 2$ is different from the case $k=1$ in important ways.  First, only for $k=1$ can the density $f_{1}(y)$ of $1\text{-}\max(X) \equiv \max (X)$ be factored into a product of the form $\phi(y) G_1(y)$ where $\phi(y)$ is the standard Guassian density and $G_1(y)$ is nondecreasing.   This factorization was a key innovation in \cite{chern:chet:kato:anti} and it follows that $f_{1}(y) = \phi(y) G_1(y) = \phi(y) G_1(y) \( \int_{y}^\infty \phi(t) dt \)^{-1}\int_{y}^\infty \phi(t) dt \leq \phi(y)\( \int_{y}^\infty \phi(t) dt\)^{-1}  \int_{y}^{\infty} \phi(t) G_1(t)dt = M(y) \Pr( \max(X) \geq y)$ where $M(y)$ is the univariate Mills ratio.  Thus, the univariate Mills ratio and  bounds for probabilities of large deviations of $\max(X)$ can be used to construct anticoncentration bounds.  Second, only for $k=1$ is $\kmax(X)$ is convex in $X$, ruling out any method for constructing desired anticoncentration bounds for $k\geq 2$ which would have needed to rely on convexity of $\kmax(X)$.

\

\noindent One potential outline for obtaining anticoncentration bounds for $k\geq 2$ which would indeed depend on $p$ is to appeal to Nazarov's anticoncentration bounds from \cite{nazarov:anticonc}.  To do so, note that there exists a Gaussian random vector $W \in \mathbb R^{_pC_k}$ with components denoted by $W_A$ and indexed by subsets $A \subseteq \{1,...,p\}$ with $|A| = k$ such that $\Ep[W_A] = \Ep[ \min_{j \in A}(X_j)]$ and $\text{cov}(W_A, W_{A'}) =  \text{cov}(\min_{j \in A}(X_j) \min_{j \in A'}(X_j))$, with the intention of using $\max(W)$ to approximate $\kmax(X)$.  The value of Nazarov's bound is that through recentering, it can be made applicable to obtain an anticoncentration bound for $\max(W)$ despite the fact that $W_A$ typically do not have mean 0.  Nazarov's bounds imply $\Pr( \max(W) \in [y,y+\eps]) \leq (\eps/\sqrt{ \min(\text{var}(W))})(\sqrt{ 2\log{_pC_k}}+2)$, where the minimum in $\min(\text{var}(W))$ spans over diagonal elements of $\text{var}(W)$.  This bound is dimension-dependent in that it depends on $p$.  Interestingly, this bound is seen to be slightly tighter in $k$ than the one in Theorem 1 after noting that $\log {_pC_k} \asymp k \log p$; though $\min(\text{var}(W))$ may be made smaller with increasing $k$.  The analysis would continue by comparing $\Pr( \max(W) \in [y, y+ \eps])$ to $\Pr( \kmax(X) \in [y , y+ \eps])$ using Gaussian comparison techniques developed in eg. \cite{chern:chet:kato:anti}.  However, this line of analysis results in dimension-dependent bounds is outside the scope of the search for dimension-free bounds.

\

\noindent  As seen in the discussion on using Nazarov's bounds in the previous paragraph, relaxing $\Ep[X_j]=0$ comes with potential loss of generality.  On the other hand, the case of more general diagonal entries in $\Sigma$, ie $\Ep[X_j^2] = \sigma_j^2$, where $\sigma_j^2$ are positive real numbers, can be reduced to the case $\Ep[X_j^2]=1$ by the same methods as were used \cite{chern:chet:kato:anti} in the case $k=1$ with a new bound depending on the minimum of the $\sigma_j^2$ but still dimension-free.

\

\noindent
The key idea here for handling the case $k\geq 2$ relative to previous literature is to compare $\kmax(X)$ to a newly defined auxiliary random variable, $k$-$\widetilde \max(X)$, constructed as follows.
For nonempty subsets $A\subseteq \{1,...,p\}$ let $\bar X_A = \frac{1}{|A|} \sum_{j\in A} X_j$.  For each $A \subseteq \{1,...,p\}$, let $\iota(A)$ be a randomly chosen, uniformly distributed, element of $A$, independent of $X$ and of $\iota(A')$ for all other $A' \subseteq \{1,...,p\}$.
Let $A^* \in \text{arg} \max_{A \subseteq \{1,...,p\}, |A|=k} \bar X_A$, with $A^*$ chosen uniformly at random from the $\text{arg} \max$ set if it is not a singleton.  
Let $\iota^*= \iota(A^*)$.  Define
$$ k\text{-}\widetilde \max(X) = X_{\iota^*}.$$

\noindent
Heuristically, $k$-$\widetilde \max(X)$ is a randomized relative of $\kmax(X)$.  Note that there is the coupling inequality $\Pr(k\text{-}\widetilde \max(X) = \kmax(X)) \geq \frac{1}{k}$. 

\

\noindent This randomization is sufficiently regularizing so that the corresponding density $\tilde f_k(y)$ can in fact be expressed in the form $\phi(y) \tilde G_k(y)$ where $\tilde G_k(y)$ is nondecreasing.  This recovers the applicability of the techniques in \cite{chern:chet:kato:anti} described above to obtain an anticoncentration bound for $\ktildemax (X)$, which subsequently translates to an anticoncentration bound for $\kmax (X).$  This factorization is stated formally in the next lemma.

\

\noindent
\textrm{LEMMA 1.} \textit{Let $X$ be as in the statement of Theorem 1.  Then $\ktildemax(X)$ is absolutely continuous with respect to Lebesgue measure and has density 
$\tilde f_k(y) = \phi(y) \tilde G_k(y)$ where $\phi(y)$ is the standard Gaussian density and
where }
$$\tilde G_k(y) =  \frac{1}{k}\sum_{j=1}^p
\Pr( j \in A^* | X_j = y).$$  \textit{Furthermore, $\tilde G_k(y)$ is nondecreasing in $y$. }

\

\noindent 
\textbf{Proof of Lemma 1}

\

\noindent
Absolute continuity follows from  $\Pr(\ktildemax(X)\in B)\leq k\Pr(X_j \in B )$ for any $j$ and Borel set $B$.  By standard reductions, following reasoning in the proof of Lemma 6 in \cite{chern:chet:kato:anti}, the next expression is well defined and provides a version of the desired density $\tilde f_k(y)$:
$$\tilde f_k(y) = \phi(y) \sum_{j =1}^p \Pr( j = \iota^* | X_j = y).$$
Next, note that $ \sum_{j=1}^p \Pr( j = \iota^* | X_j = y)=  \sum_{j=1}^p \Pr( j \in A^*, j = \iota^*|X_j = y)  =  \sum_{j=1}^p\Pr(j = \iota^*| j\in A^*, X_j=y) \Pr( j \in A^* | X_j = y) = \frac{1}{k}  \sum_{j=1}^p \Pr( j \in A^* | X_j = y)  = \tilde G_k(y)$.  Thus $\tilde f_k(y) = \phi(y) \tilde G_k(y)$ as stated in Lemma 1. 

\

\noindent
To show that $\tilde G_k(y)$ is nondecreasing, note $\Pr(j \in A^*|X_j = y) = \Pr( \min_{l\in A}(X_l) \leq X_j \text{ for all }A \subseteq \{1,...,p\}, |A| = k | X_j = y).$
Let $V_{jl} = X_l - \Ep[X_l X_j]X_j$, ie. $V_{jl}$ are the residuals from the least squares projection of $X_l$ on $X_j$ and are jointly independent of $X_j$ as $X$ is jointly Gaussian. Furthermore, for any given set $A$,
$$\min_{l\in A}(V_{jl} + \Ep[X_lX_j] y) \leq y  \ \Leftrightarrow  \ \min_{l \in A} (V_{jl} + (\Ep[X_lX_j] - 1)y) \leq 0.$$
Independence of $\{V_{jl}\}_{l=1,...,p}$ from $X_j$ and the fact that $(\Ep[X_lX_j] - 1) \leq 0$ then implies $\Pr( j \in A^* | X_j = y)$ is nondecreasing in $y$.  \ \QED

\

\noindent
\textbf{Proof of Theorem 1} 

\

\noindent
Given the conclusion of Lemma 1, then 
$$\tilde f_k(y) = \phi(y) \tilde G_k(y) = \phi(y) \tilde G_k(y) \( \int_y^\infty \phi(t) dt \)^{-1}\( \int_y^\infty \phi(t) dt \)$$
$$ \leq  \phi(y) \( \int_y^\infty \phi(t) dt \)^{-1}\( \int_y^\infty \phi(t) \tilde G_k(t) dt \) = M(y) \Pr( \ktildemax(X) \geq y).$$ Because $\ktildemax(X) \leq \max(X)$ it holds that $\Pr( \ktildemax(X) \geq y) \leq \Pr( \max(X) \geq y)$ and therefore $$ \tilde f_k(y) \leq M(y) \Pr( \max(X) \geq y).$$ Next, under the condition that $\Ep[X_j^2]=1$ for all $j$, \cite{chern:chet:kato:anti} prove that $\Pr( \max(X) \in [y,y+\eps]) \leq 2\eps(1+ \Ep[\| X \|_\infty])$ using only the property that the density, $f_1(y)$, of $\max(X)$, satisfies $f_1(y) \leq M(y) \Pr( \max(X) \geq y)$.  This same property holds for $\tilde f_k(y)$ and thus also $\Pr(\ktildemax(X) \in [y,y+\eps]) \leq 2\eps(1 + \Ep [ \| X \|_\infty] ).$  Finally, by construction, 
$\Pr(\kmax(X) \in [y,y+\eps]) \leq k \Pr(\ktildemax(X) \in [y,y+\eps])$, and Theorem 1 follows. \ \QED

\

\noindent
\textbf{Discussion of applications of Theorem 1 to multiple testing problems}  

\

\noindent As an application of Theorem 1,  
consider the problem of testing a collection of $p$ statistical hypotheses $H_{01},...,H_{0p}.$  When $p \geq 1$, conducting separate testing procedures for each $H_{0j}$ at level $\alpha \in (0,1)$ independently may lead to the classical multiple testing problem that the family-wise error rate exceeds $\alpha$, ie. $\mathsf{FWER} = \Pr( H_{0j} \text{ rejected for some true null } j ) > \alpha$, even if $\Pr( H_{0j} \text{ rejected }) \leq \alpha$ individually for all $j$ corresponding to a true $H_{0j}$.  To address this, researchers have designed many joint testing procedures aimed at controlling family-wise error rate. 

\

\noindent
In many applications, $\mathsf{FWER}$ is too stringent a notion of error rate control.  There are several alternatives to \textsf{FWER} which allow for control in the tradeoff between power and tolerance for false positives in multiple testing problems with many hypotheses.  These include $k$-family-wise error rate, $k\text{-}\mathsf{FWER}$ $= \Pr( H_{0j} \text{ rejected for at most $k$ true nulls } j)$.  $k$-\textsf{FWER} is a natural starting point for discussion due to its simplicity, and procedures controlling other error rate notations can sometimes be built $k$-\textsf{FWER}-controlling procedures; eg. the \textsf{FDP}-controlling procedure of \cite{romano:wolf:2007}.

\

\noindent
One common method for controlling $k$-\textsf{FWER} is a step-down-based method in Algorithm 2.1 in \cite{romano:wolf:2007}. Their procedure depends only on test statistics $T_1,...,T_p$  corresponding to $H_{01},...,H_{0p}$ and estimated critical values $\hat c_{K}(1-\alpha,k) \text{ for each } K \subseteq \{1,...,p \}, |K| \geq k$, which preferably estimate upper bounds over $1-\alpha$ quantiles of $\kmax(T_j: j \in K)$ under distributions in the intersection null, $\cap_{j \in K} H_{0j}$.  Note, $\kmax(T_j: j \in K)$ refers to $\kmax$ applied to the vector with components $T_j$ with $j\in K$.   Let $I\subseteq \{1,...,p\}$ be the set of true nulls.  Theorem 2.1 of \cite{romano:wolf:2007} proves that under an additional monotonicity condition that $\hat c_{K}(1-\alpha,k) \geq \hat c_I (1-\alpha,k)$ for $K \supseteq I$, their Algorithm 2.1 results in $k$-\textsf{FWER} $\leq \Pr(\kmax(T_j: j \in I) > \hat c_{I}(1-\alpha,k))$.

\

\noindent
Because large-$p$ settings are precisely those for which control of generalized error rates is highly relevant, a natural question is: to what extent do $\hat c_{K}(1-\alpha,k)$ result in favorable statistical properties for $k$-\textsf{FWER} as a function of both $n$ and $p$ when based on bootstrap estimates of order statistics of vectors with components being subsets of $T_1,...,T_p$?   Note,
\cite{romano:wolf:2007} discussed conditions in which $k$-\text{FWER} is asymptotically controlled with bootstrap-based $\hat c_{K}(1-\alpha,k))$ in a frame in which $n\rightarrow \infty$ while $p$ was fixed, thus suppressing dependence of $k$-\textsf{FWER} on $p$ from notation
For illustration, suppose $U_1,...,U_n\in \mathbb R^p$ is a collection of $n$ independent, identically distributed random Gaussian vectors with unknown mean $\mu \in \mathbb R^p$.  Write $U_{ij}$ and $\mu_j$ for the $j$th component of $U_i$ and $\mu$. Suppose $H_{0j}$ are given by $H_{0j} : \mu_j \leq 0$ versus $\mu_j > 0$.  Suppose covariance of the $U_i$ is unknown with the exception, for the sake of simplicity, that the diagonal elements of the covariances are finite and known (and assumed =1). 
Let $T_{j} = n^{-1/2} \sum_{i=1}^n U_{ij}$. Recall that the test that rejects for large $T_j$ is uniformly most powerful for testing $H_{0j}$.  Let $U^*_1,...,U_n^*$ be an empirical bootstrap sample from $U_1,...,U_n$, and let $T^*_j = n^{-1/2} \sum_{i=1}^n (U_{ij}^* -n^{-1}\sum_{q=1}^nU_{qj}).$  Let $\hat c_K(1-\alpha, k)$ be the $1-\alpha$ quantile of $\kmax(T^*_j: j \in K)$ (conditioned on $U_1,...,U_n$, calculated over draws of $*$).  Apply Algorithm 2.1 of \cite{romano:wolf:2007} to yield decisions $D_j \in \{Reject, \ Fail \ to \ Reject\}$ for each $H_{0j}$.

\

\noindent
Properties of $D_j$ follow from analyzing $\Pr( \kmax(T_j:j \in I) > \hat c_I(1-\alpha,k))$.   
Let $\beta =  \left ( q_{1-\alpha}(\kmax(T_j : j \in I ))-\hat c_I(1-\alpha,k) \right )$, and let $\gamma, \delta \geq 0 $ be such that $\Pr( \beta \geq \gamma ) \leq\delta$.   By application of Theorem 1 above, $\Pr( \kmax(T) \geq \hat c_I(1-\alpha,k) ) \leq \alpha + 2 k \gamma (1 + \Ep[ \|U\|_\infty]) + \delta.$
Application of Theorem 2.2 of \cite{romano:wolf:2007}, along with the fact that their monotonicity condition is satisfied in the above construction, gives
$$k\text{-}\mathsf{FWER} \leq \alpha + 2k\gamma(1+ \Ep[\|U\|_\infty]) + \delta.$$

\noindent
Due to its use in deriving the above expression, a dimension-free bound for the anticoncentration of Gaussian order statistics can improve finite sample upper bounds on deviations of $k\text{-}\mathsf{FWER}$ from $\alpha$.  The dimension-free property becomes important when there is enough correlation in the components of $U_i$.   Although an exhaustive treatment of conditions which lead to control of $\gamma$ and $\delta$ are outside of the current scope, note that following arguments like those in \cite{clt:empboot:highdim} can then give more explicit bounds on $\gamma, \delta$.  Note that $\Pr( \beta \geq \gamma ) \leq\delta$ is a concentration property rather than an anticoncentration property, and is implied by $\Pr( \sup_{y\in \mathbb R}| \Pr( \kmax(T_j:j \in I) < y ) - \Pr(\kmax( T_j^*:j\in I )<y |U_1,...,U_n) | < \gamma |) \geq 1- \delta$.   Note also, the above bound for $\kmax$ is applied in the same way that \cite{chet:chern:kato:multboot} studied bootstrap estimates $\hat c_K(1- \alpha, 1)$ in the case $k=1$ leading to bootstrap-based procedures controlling \textsf{FWER}.     

\

\noindent
\textbf{Additional remarks}  

\

\noindent There are several potential avenues for future research expanding on the bounds of Theorem 1.  First, understanding anticoncentration properties of the $k$th largest local maximum of almost-surely smooth Gaussian processes $X$ with components indexed on a disjoint union of $p$ intervals $\bigsqcup_{j=1}^p [0,1]_j$ can lead to bootstrap uniform confidence bands which cover an unknown function $g:\bigsqcup_{j=1}^p [0,1]_j \hspace{-1mm}\rightarrow \mathbb R$ on all but $k$ intervals $[0,1]_j$ with probability at least $1-\alpha$ using methods like those developed in \cite{Chernozhukov_2014}.  Second, the randomization technique leading to the definition of $\ktildemax(X)$ is distinct from other randomization or symmetrization techniques used in the study of empirical processes (eg. Lemma 2.3.1 in \cite{vdV-W}), and may have additional applications outside of Theorem 1.  
Finally, the dimension-free anticoncentration bounds can potentially lead to additional improvements in understanding $k$-\textsf{FWER} control with a bootstrap-based procedure for non-Gaussian data if coupled with interpolation and Gaussian comparison techniques like those developed in \cite{CHERNOZHUKOV20163632}, \cite{chernozhukov2019improved} and \cite{chernozhukov2021nearly}.   For example, if $U_i$ themselves are not jointly Gaussian, then interpolation techniques can be used to compare the distributions of $\kmax(Q)$ to $\kmax(T)$, where $Q$ is defined to be a Gaussian random vector with components indexed by $j=1,...,p$ with $\Ep[Q_j] = 0$ and $\Ep[Q_jQ_l]=\Ep[T_{j} T_{l}].$


\bibliographystyle{apalike}
\bibliography{dkbib1}

\end{document}